\documentclass[11pt,a4paper]{article}

\usepackage{amssymb,graphicx,epsfig,color}

\usepackage[a4paper, total={6in, 8in}]{geometry}
\usepackage{amsmath,amssymb,graphicx,epsfig,color,amsthm}
\usepackage{hyperref}

\renewcommand*{\d}{\mathrm{d}}

\newtheorem{Theorem}{Theorem}
\newtheorem*{Theorem*}{Theorem}

	\theoremstyle{remark}
\title{{\bf Spirals of Riemann's Zeta-Function}\\ {\large --- Curvature, Denseness, and Universality ---}}

\author{Athanasios Sourmelidis \& J\"orn Steuding}
\begin{document}

\date{}

\maketitle

\begin{abstract}
This article deals with applications of Voronin's universality theorem for the Riemann zeta-function $\zeta$. 
Among other results we prove that every plane smooth curve appears up to a small error in the curve generated by the values $\zeta(\sigma+it)$ for real $t$ where $\sigma\in(1/2,1)$ is fixed. 
In this sense, the values of the zeta-function on any such vertical line provides an atlas for plane curves.
In the same framework, we study the curvature of curves  generated from $\zeta(\sigma+it)$ when $\sigma>1/2$ and we show that there is a connection with the zeros of $\zeta'(\sigma+it)$. 
Moreover, we clarify under which conditions the real and the imaginary part of the zeta-function are jointly universal. 
\end{abstract}
{\small{\noindent {\sc Keywords:} Riemann zeta-function, Voronin's universality theorem, curvature, Jensen's function }\\
{\sc Funding:} AS was supported by the Austrian Science Fund (FWF) project M 3246-N.}
\section{Introduction}

About fifty years ago, Sergey Voronin \cite{voronin}\footnote{Note that the paper in question was submitted in 1974 and only published in 1975.} proved his celebrated universality theorem for the Riemann zeta-function $\zeta(s)$, $s=\sigma+it$ being a complex variable. This astonishing result states that, roughly speaking, certain shifts of the zeta-function approximate every zero-free analytic function, defined on a sufficiently small disk (see also \cite{karatsuba}). In this note we discuss a few new consequences of this remarkable property with respect to the curves given by the values of $\zeta(\sigma+it)$ as $\sigma$ is fixed and $t$ ranges through the set of real numbers or some subinterval. These curves look like {\it spirals} when $t$ is from a bounded range, and we will use that word from time to time in what follows. The implications of universality that we consider here are by no means complicated in nature though they seem to have been overlooked so far. 
\smallskip

For our purpose we recall the universality theorem \cite{voronin} in a stronger form: {\it Suppose that ${\mathcal K}$ is a compact subset of the strip $1/2<{\rm{Re}}\,s<1$ with connected complement, and let $g(s)$ be a non-vanishing continuous function on ${\mathcal K}$ which is analytic in the interior of ${\mathcal K}$. Then, for every $\epsilon>0$,}
$$
\liminf_{T\to\infty}{1\over T}{\rm{meas}}\left\{\tau\in[0,T]\,:\,\max_{s\in{\mathcal K}}\vert \zeta(s+i\tau)-g(s)\vert<\epsilon\right\}>0
$$
(see \cite{steu}). The main differences to Voronin's original statement in \cite{voronin} are the positive lower density of the set of shifts $\tau$ (which is already implicit in Voronin's proof) and the rather general set ${\mathcal K}$ where Voronin considered only disks; this is first apparent in Steve Gonek's thesis \cite{gone} and later in Bhaskar Bagchi's thesis \cite{bagchi}. The topological restriction for ${\mathcal K}$ follows from Mergelyan's approximation theorem and its limitations (see \cite{mergelyan} and \cite{steu}, p. 107). We will also make use of the following observation due to Johan Andersson \cite{ander}: {\it If $\mathcal{K}$ has empty interior, then the target function $g$ in the universality theorem is allowed to have zeros}.

Assuming the Riemann Hypothesis, Gonek \& Hugh Montgomery \cite{gonek} showed that the parametrized curve $t\mapsto \zeta(1/2+it)$ turns in the clockwise direction for all sufficiently large $t$ or, in other words, the spiral has negative curvature (see Figure 1 for an illustration).\footnote{In this context it is interesting to notice that Daniel Shanks \cite{shanks} conjectured that $t\mapsto \zeta(1/2+it)$ {\it approaches the origin (at the nontrivial zeros) mostly from the third quadrant}; this was proved by Akio Fujii \cite{fujii}. The conditional curvature result of Gonek \& Montgomery matches this scenario, in particular in combination with Harold Edwards' observation that the values $\zeta(1/2+it)$ {\it lie most of the time in the right half-plane} \cite{edwards}, which was verified (in some sense) by Justas Kalpokas \& the second author \cite{kalpokas}.} A simple application of the universality theorem shows that the behaviour to the right of the critical line can be rather different.  
\bigskip

\begin{figure}[h]
\includegraphics[height=3.72cm]{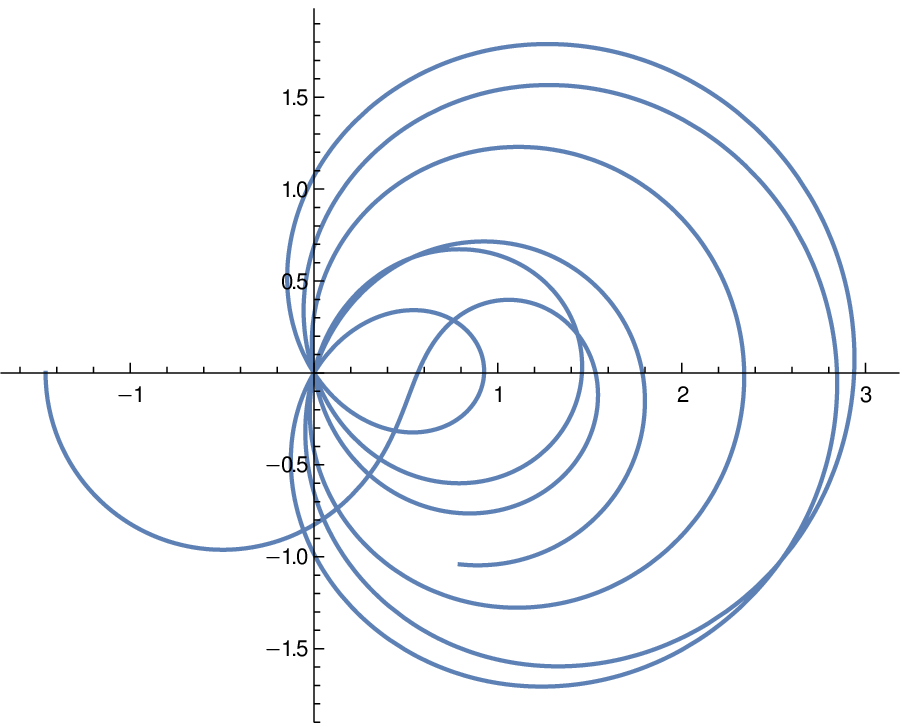}\qquad \includegraphics[height=3.72cm]{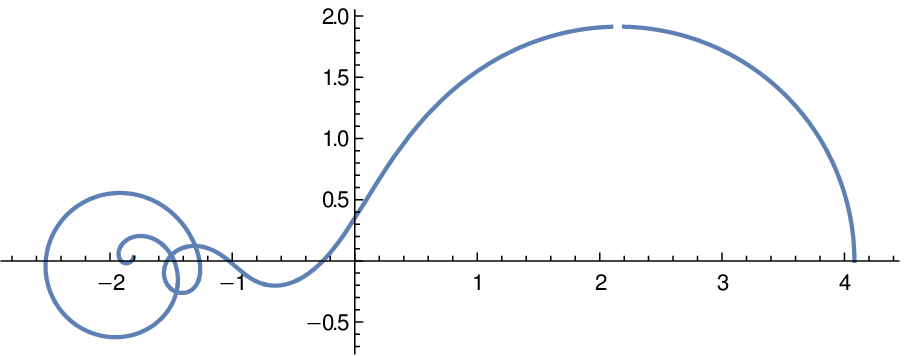}
\caption{Left: $\zeta(1/2+it)$ for $0\leq t\leq 40$; right: $\zeta''/\zeta'(1/2+it)$ for the same range of $t$. The negative real part of $\zeta''/\zeta'(1/2+it)$ for $t\geq t_0=2.75\ldots$ corresponds to the clockwise direction of the spiral on the left.}
\end{figure}
\bigskip

For this purpose, we consider any parametrized curve 
$$
{\mathcal C}={\mathcal C}(\sigma,{\mathcal I})\ :\qquad I\ni t\mapsto \zeta(\sigma+it),
$$
where ${\mathcal I}=[a,b]$ is an interval, $\sigma\in(1/2,1)$ is fixed. 
Then ${\mathcal K}=\{\sigma+it\,:\, t\in {\mathcal I}\}$ is a compact set in $1/2<{\rm{Re}}\,s<1$ with connected complement and empty interior. Hence, applying the universality theorem with $g(s)=\zeta(s)$, implies the reappearance of the curve ${\mathcal C}$ within $\zeta(\sigma+i\mathbb{R})$ up to an invisible error of size $\epsilon$ on a set of positive lower density.    

However, we may also consider the {\it inverse} curve
$$
{\mathcal C}^*={\mathcal C}^*(\sigma,{\mathcal I})\ :\qquad {\mathcal I}\ni t\mapsto \zeta(\sigma+i(a+b-t)).
$$
This curve consists of the same points as ${\mathcal C}$, only its direction is inverted. Thus, applying the universality theorem with $g(\sigma+it)=\zeta(\sigma+i(a+b-t))$ now yields the appearance of ${\mathcal C}^*$ within $\zeta(\sigma+i\mathbb{R})$ up to an invisible error of size $\epsilon$ for a set of positive lower density. 

We observe that the curvature of ${\mathcal C}^*$ is the negative of the curvature of ${\mathcal C}$. So the curvature of the parametrized curves $t\mapsto \zeta(\sigma+it)$ for $\sigma\in(1/2,1)$ takes both signs infinitely often as $t\to\infty$, and both signs appear for a set of positive lower density. This is completely different to the situation on the critical line (if the Riemann Hypothesis is true). 

Note that the curvature changes its sign when $t$ tends to $-\infty$ or $\zeta(\sigma+it)$ is replaced by its conjugate $\zeta(\sigma-it)$. But this change of orientation is of different nature as the reappearance of ${\mathcal C}^*$ in $\zeta(\sigma+i\mathbb{R})$ (up to a small error). The phenomenon in the background is, probably, an almost periodicity property of $\zeta(\sigma+it)$ for $\sigma>1/2$ in combination with universality. Computer experiments occasionally show an instance of a positive curvature in these spirals (see Figure 2).

\begin{figure}[h]
\includegraphics[height=6.8cm]{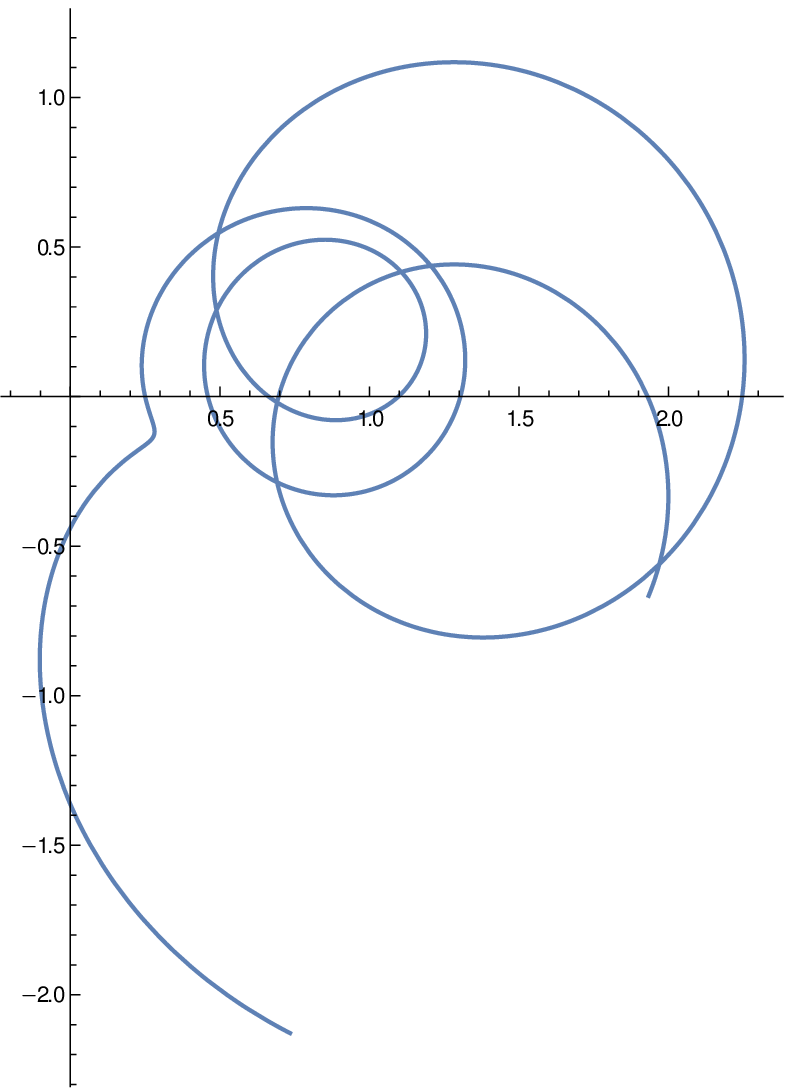}\qquad\includegraphics[height=6.8cm]{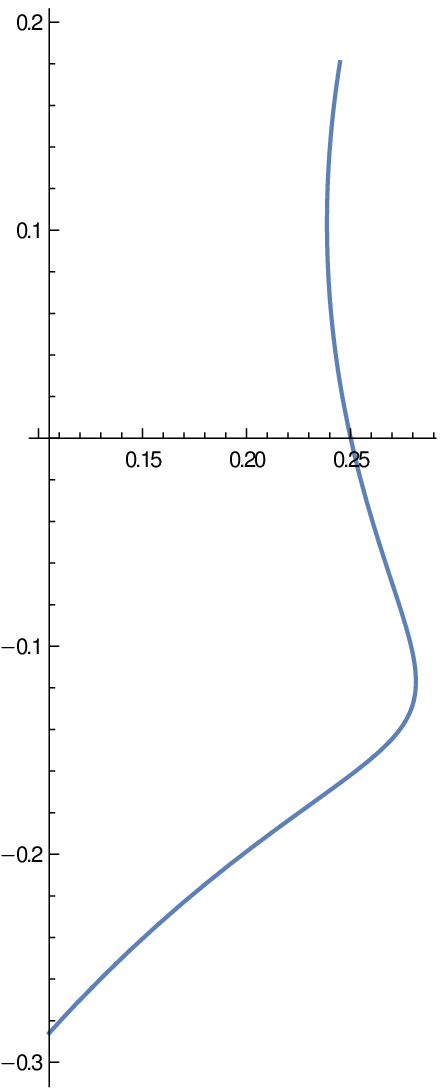}\qquad\includegraphics[height=6.8cm]{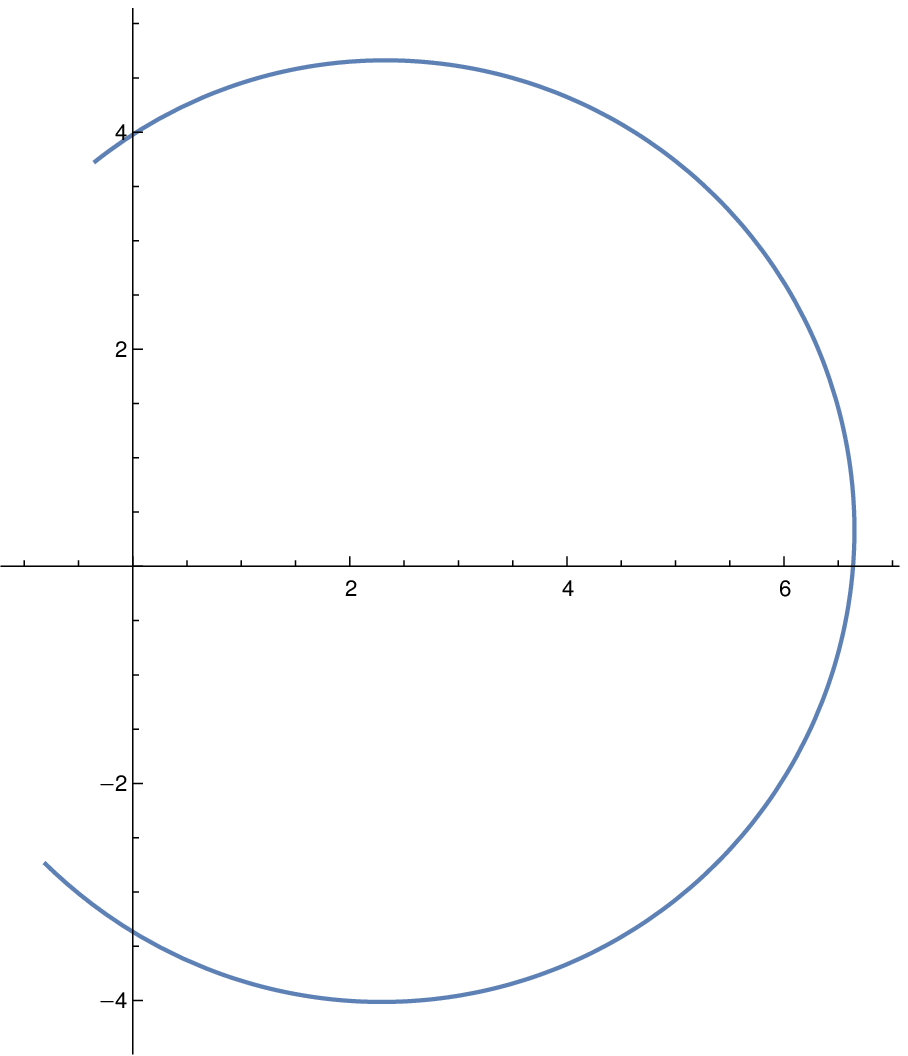}
\caption{Left: $\zeta(3/4+it)$ for $110\leq t\leq 120$; middle: the same for the range $111.2\leq t\leq 111.7$; right: $\zeta''/\zeta'(3/4+it)$ for the same range of $111.2\leq t\leq 111.7$, where a positive real part indicates a positive curvature.}
\end{figure}

\section{Statement of the Main Results}

In the sequel we assume every curve to be {\it smooth}, by which we mean that the plane and space curves under investigation have a parametrization with at least second and third order continuous derivatives, respectively. We also assume that they are {\it regular} which means that their first derivative is non vanishing. We begin with a remarkable generalization of the result from the introduction.

\begin{Theorem}\label{1}
Let $\sigma\in(1/2,1)$ and $\epsilon>0$ be fixed. Then the following statements are true:
\begin{enumerate}
\item  Every plane curve is up to an error of size $O(\epsilon)$ and affine translation contained in the graph of the curve $\mathbb{R}\ni t\mapsto\zeta(\sigma+it)\in\mathbb{C}$. 

\item  Every space curve is encoded in the curve $\mathbb{R}\ni t\mapsto\zeta(\sigma+it)\in\mathbb{C}$ up to an error of size $O(\epsilon)$.
\end{enumerate}
 \end{Theorem}

\noindent Here, of course, we consider a plane curve in the euclidean plane via $\mathbb{R}^2\simeq \mathbb{C}$ also as a curve in the complex plane. 

By the first statement of the theorem, the values of the zeta-function on any vertical line in the right open half of the critical strip provide an atlas for plane curves (similar as a normal number contains any finite pattern of digits).\footnote{A paper by Elias Wegert \& Gunter Semmler \cite{wegert} contains another interesting though different application of universality to plane curves.} A short note with this result and a sketch of its proof has been submitted \cite{19}.
\bigskip

\begin{figure}[h]
\includegraphics[height=8cm]{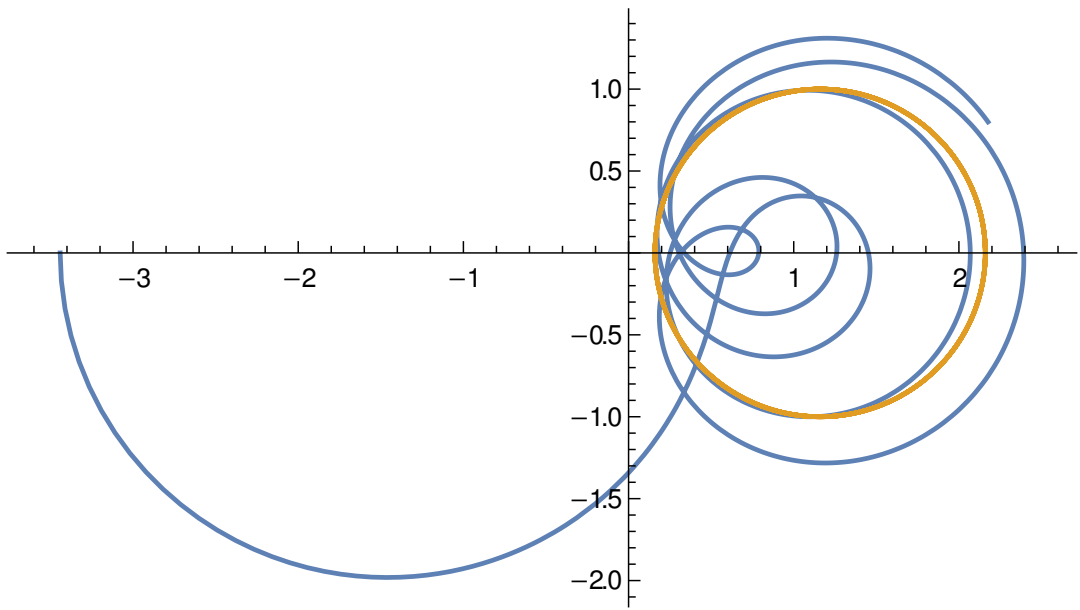}
\caption{The spiral $\zeta(3/4+it)$ for $0\leq t\leq 35$ containing an approximation of a circle (in yellow).}
\end{figure}
The second statement relies on a simultaneous approximation of the curvature and the torsion of a space curve. The proof (in Section \ref{11}) is indeed an instance of a joint universality phenomenon for the real and the imaginary part of the zeta-function. 
Interestingly, the question whether the real and the imaginary part of the zeta-function can simultaneously approximate admissible target functions depends mainly on the set of definition of the target functions and can be answered positively and negatively:

\begin{Theorem}\label{2}
The real and imaginary part of the zeta-function are jointly universal with respect to approximation of continuous real-valued functions defined on a compact set ${\mathcal K}\subseteq\lbrace s\in\mathbb{C}:1/2<{\rm Re}s<1\rbrace$ with connected complement if, and only if, $\mathcal{K}$ has empty interior. More precisely, assume that $\mathcal{K}$ has empty interior and $f, g\,:\,{\mathcal K}\to\mathbb{R}$ are continuous.	Then, for every sufficiently large $T$ and every $\epsilon>0$, there exists $\tau\in[T,2T]$ such that
$$
\left\{\begin{array}{l}\max_{s\in{\mathcal K}}\vert{\rm{Re}}\,\zeta(s+i\tau)-f(s)\vert<\epsilon,\\ \max_{s\in{\mathcal K}}\vert{\rm{Im}}\,\zeta(s+i\tau)-g(s)\vert<\epsilon;\end{array}\right.
$$
moreover, the set of such $\tau$ has positive lower density. On the other hand, even if we assume any order of differentiability for $f$ and $g$, the above statement can not be true if $\mathcal{K}$ contains an open disk.
\end{Theorem}

The positive statement of the above theorem shows that also curiosities as the Peano curve or other space-filling curves can be approximated as long as they have a continuous representation  (see \cite{sagan}).

Motivated by the aforementioned result of Gonek \& Montgomery about the clockwise direction of the spiral associated with the critical line we also consider the curves resulting from the values of the zeta-function on vertical line segments to the left of the critical line and to the right of the critical strip, respectively. It seems that the behaviour in the right open half of the critical strip is exceptional:

\begin{Theorem}\label{3} We consider the parametrized curve $t\mapsto \zeta(\sigma+it)$ for fixed $\sigma$ and any sufficiently large $t\geq t_0(\sigma)>0$.
\begin{enumerate}
\item There is a constant $A\in[3,4]$ such that the curvature is negative for $\sigma\geq A$.

\item The curvature is negative for $\sigma\leq 0$ and under Riemann hypothesis also for $\sigma\leq 1/2$.

\item There exists $\theta>0$ such that the curvature has infinitely many sign changes for $\sigma\in(1/2,1+\theta)$.
\end{enumerate}
\end{Theorem}

{ We will see below that the curvature of the curve $t\mapsto \zeta(\sigma+it)$ has the same sign as
$$
\mathrm{Re} {\zeta''\over \zeta'}(\sigma+it).
$$
Hence, to prove Theorem \ref{3} it suffices to examine which sign does ${\rm{Re}}\, \zeta''/\zeta'(\sigma+it)$ take along vertical lines for various values of $\sigma\neq1/2$.	
Although the last statement of Theorem \ref{3} says that the curvature can be positive an infinitude of times for any $\sigma\in(1/2,1+\theta)$, we will show that on average it is always negative.
Indeed, we have the following
}
\begin{Theorem}\label{3.5}
For every $\sigma>1/2$, the limit 
$$
\lim_{T\to\infty}\frac{1}{T}\int_0^T\mathrm{Re} {\zeta''\over \zeta'}(\sigma+it)
$$ 
exists and it represents an increasing and eventually constant function of $\sigma$ with maximum value $-\log 2$.
\end{Theorem}

{As a matter of fact, we will prove that the aforementioned limit at a point $\sigma$ is the derivative of the so-called {\it Jensen function} of $\zeta'$, showing therefore a connection with the relative frequency of zeros of $\zeta'$ in the right of the vertical line $1/2+i\mathbb{R}$.}

Lastly, we consider a result of Ram\={u}nas Garunk\v{s}tis \& the second author \cite{garunk} who showed, under assumption of the Riemann Hypothesis, that $\zeta(\sigma+i\mathbb{R})$ is not dense in $\mathbb{C}$ for any fixed $\sigma<1/2$. Michel Lapidus \cite{lapi} proved that also the converse is true; see also the recent book of Hafedh Herichi \& Lapidus \cite{lapidus}. We can strengthen these results slightly as follows:     

\begin{Theorem}\label{4}
The Riemann Hypothesis is true if, and only if, $\zeta(\sigma+i\mathbb{R})$ is nowhere dense in $\mathbb{C}$ for any fixed $\sigma<1/2$.
\end{Theorem}

\noindent Harald Bohr \& Richard Courant \cite{bohr} proved that $\zeta(\sigma+i\mathbb{R})$ is dense in $\mathbb{C}$ for fixed $\sigma\in(1/2,1]$. Of course, this result also follows from universality theorem (by choosing a constant target function). For the critical line, however, it is unkown whether $\zeta(1/2+i\mathbb{R})$ dense in the complex plane or not. Universality does not apply to the critical line because of too many zeta zeros.
\medskip

In the following four sections we give the proofs of these results. All our reasonings apply to more general zeta- and $L$-functions. What is necessary is a suitable universality theorem. Concerning the behaviour of the curvature and its sign changes to the right of the critical line, we conclude in Section 4 with the determination of the so-called Jensen function of $\zeta'$ (plus an appendix).  

\section{Proof of Theorem \ref{1}}\label{11}

It is a well known fact that a smooth plane curve is determined by its curvature. In fact, the fundamental theorem of the local theory of curves states that {\it given differentiable functions $\mathfrak{t},\kappa$ defined on an interval ${\mathcal I}$ satisfying $\kappa(t)>0$ for $t\in{\mathcal I}$, there exists a regular parametrized curve $\gamma\,:\,{\mathcal I}\to \mathbb{R}^3$ such that $t$ is the arclength (i.e. the euclidean norm of $\gamma'(t)$ equals $1$ for every $t\in\mathcal{I}$), $\kappa(t)$ is the curvature, and $\mathfrak{t}(t)$ is the torsion of $\gamma$. Moreover, any other curve $\tilde{\gamma}$, satisfying the same conditions, differes from $\gamma$ by a rigid motion} (meaning that there exists an orthogonal linear map $\ell$ and a vector $v$ such that $\tilde{\gamma}=\ell\circ \gamma+v$); see, for example, \cite{carmo}. For plane curves there is no torsion, hence this result implies that a plane curve $\gamma$ is indeed determined by its curvature; in this case, the curvature is not restricted to be positive. The proof relies on solving a system of certain differential equations, the so-called {\it Frenet's equations}, and this plays also a certain role in the form of the integral equations  (\ref{111})  and \eqref{victor2} below. 

We return to the first statement of Theorem \ref{1}. Let $\kappa$ be the curvature of a plane curve ${\mathcal C}$. Define
\begin{equation}\label{111}
\vartheta(u):=\int_0^u\kappa(t)\d t.
\end{equation}  
Then, a model of the curve ${\mathcal C}$ with curvature $\kappa$ and arclength $t$ in the complex plane is given by the parametrization
$$
t\mapsto g(t):=\int_0^t\exp\big(i\vartheta(u)\big)\d u,
$$ 
where $t$ ranges through some interval ${\mathcal I}$. By the universality theorem, for every $\epsilon>0$, there exists $\tau>0$ such that
\begin{equation}\label{24}
\max_{t\in{\mathcal I}}\left\vert \zeta(\sigma+it+i\tau)-g(t)\right\vert<\epsilon.
\end{equation} 
In view of the positive lower density for the real shifts $\tau>0$ that lead to the desired approximation of the target function it follows that any plane curve appears infinitely often, up to a tiny error, in any curve $\zeta(\sigma+i\mathbb{R})$ with any fixed $\sigma\in(1/2,1)$ (even with positive lower density). In this sense, {\it the zeta-function provides a single plane curve that contains all plane curves with an error too small to be seen with the naked eye!}\,\footnote{The Planck length is about $1.6\cdot 10^{-36}$ meters and, according to quantum mechanics, one cannot {\it see} anything smaller than this tiny quantity.}
\smallskip

Next we consider space curves. If $c\,:\,{\mathcal I}\to\mathbb{R}^3$ is a parametrized curve with respect to its arclength $t$ of curvature $\kappa(t)>0$ and torsion $\mathfrak{t}(t)$, then the normal vector $n(t)$ and binormal vector $b(t)$ are given by
\begin{equation}\label{victor1}
n(t)={c''(t)\over \kappa(t)}\qquad\mbox{and}\qquad b(t)=c'(t)\times n(t),
\end{equation}
respectively (where $\times$ is the vector product). In this case the curve, defined by 
\begin{equation}\label{victor2}
\tilde{c}(t)=-\int_{0}^t\frac{n'(u)+\mathfrak{t}(u)b(u)}{\kappa(u)}\d u,
\end{equation}
is congruent to $c$ up to some vector in $\mathbb{R}^3$, because ${\tilde{c}}'(t)=c'(t)$. Applying the universality theorem to the target function $t\mapsto \kappa(t)+i\mathfrak{t}(t)$ yields the desired approximation, where the quantities $\kappa$ and $\mathfrak{t}$ have in view of (\ref{victor1}) and (\ref{victor2}) to be replaced by ${\rm{Re}}\,\zeta(\sigma+it+i\tau)$ and ${\rm{Im}}\,\zeta(\sigma+it+i\tau)$.

\section{Proof of Theorem \ref{2}}

We begin with the case of continuous real-valued functions $f,g$ defined on a compact set with connected complement ${\mathcal K}\subseteq\lbrace s\in\mathbb{C}:1/2<\sigma<1\rbrace$ that has empty interior. For every sufficiently large $T$ and every $\epsilon>0$, by a direct application of the universality theorem, there exists $\tau\in[T,2T]$ such that
$$
\max_{s\in{\mathcal K}}\vert \zeta(s+i\tau)-(f(s)+ig(s))\vert<\epsilon.
$$
Separating real and imaginary part, we thus have 
\begin{align}\label{jointuniv}
\left\{\begin{array}{l}\max_{s\in{\mathcal K}}\vert{\rm{Re}}\,\zeta(s+i\tau)-f(s)\vert<\epsilon,\\ \max_{s\in{\mathcal K}}\vert{\rm{Im}}\,\zeta(s+i\tau)-g(s)\vert<\epsilon.\end{array}\right.
\end{align}
Obviously, the set of such $\tau$ has positive lower density. This proves the first assertion of Theorem \ref{2}.
\smallskip

Next we consider the case where ${\mathcal K}$ has a non-empty interior or, equivalently, that it contains an open disk $D$.
Assume that the real and imaginary part of the zeta-function do approximate two functions $f, g\,:\,{\mathcal K}\to\mathbb{R}$ in the sense of \eqref{jointuniv}. Then we can find a sequence of real numbers $\tau_n>0$ with
$$
\left\{\begin{array}{l}\sup_{s\in D}\vert{\rm{Re}}\,\zeta(s+i\tau_n)-f(s)\vert<n^{-1},\\ \sup_{s\in D}\vert{\rm{Im}}\,\zeta(s+i\tau_n)-g(s)\vert<n^{-1}.\end{array}\right.
$$
The above implies that the sequence of analytic functions $\lbrace\zeta(s+i\tau_n)\rbrace_{n\geq1}$ converges uniformly in $D$ to the function $F(s):=f(s)+ig(s)$. Hence, $F(s)$ is analytic in $D$ and,  therefore,  the partial derivatives of $f$ and $g$ as functions of two variables $(\sigma,t)$ exist and satisfy the Cauchy--Riemann equations in $D$. But this can be false in general if we choose a pair of functions $f,g$ for which the Cauchy--Riemann equations fail.

\section{Proof of Theorem \ref{3}} 

Recall the aformenentioned result of Gonek \& Montgomery \cite{gonek} who showed that the parametrized curve $t\mapsto \zeta(1/2+it)$ turns in the clockwise direction for all sufficiently large $t$. This follows from an older result due to Cem Y\i ld\i r\i m \cite{yildirim} implying a negative curvature: 
$$
\kappa(t)={{\rm{Re}}\, {\zeta''\over\zeta'}(1/2+it)\over \vert \zeta'(1/2+it)\vert}<0\qquad\mbox{for}\quad t\geq 2.76.
$$
Both results are conditional subject to the truth of the Riemann Hypothesis. 

More generally, it can be seen, similarly as in \cite{gonek}, that the curvature of the parametrized curve $t\mapsto \zeta(\sigma+it)$ is given by 
$$
\kappa(t)={{\rm{Re}}\, {\zeta''\over \zeta'}(\sigma+it)\over \vert\zeta'(\sigma+it)\vert},
$$
whenever $\zeta'(\sigma+it)\neq 0$. 
{Thus, the sign of the curvature matches the one of $\mathrm{Re}\frac{\zeta''}{\zeta'} (\sigma+it)$.}

From straightforward computations of the Dirichlet series expansion of $\zeta^{(k)}(s)$, $k=1,2$, it follows that there is a constant $A\in[3,4]$ with
\begin{equation}\label{eleminequ}
\left|1-\frac{(-1)^k\zeta^{(k)}(s)2^s}{(\log 2)^k}\right|\leq\frac{\sqrt{2}}{2},\qquad\sigma\geq A,\quad t\in\mathbb{R},\quad k=1,2.
\end{equation}
Since
$$
{\rm{Re}}\, {\zeta''\over \zeta'}(s)=-\log 2\, \mathrm{Re}\frac{\frac{\zeta''(s)2^s}{(\log 2)^2}}{\frac{-\zeta'(s)2^s}{\log 2}},
$$
we conclude from \eqref{eleminequ} that the curvature is negative for all $\sigma\geq A$ and the corresponding spirals turn in clockwise direction. 

In the half-plane $\sigma\leq0$ we can argue as Yildirim \cite{yildirim} did for $0\leq\sigma\leq1/2$. We start with the partial fraction representation 
\begin{align}\label{pfdf}
\frac{\zeta''}{\zeta'}(s)=\frac{\zeta''}{\zeta'}(0)-2-\frac{2}{s-1}+\sum_{n\geq1}\left(\frac{1}{s+a_n}-\frac{1}{a_n}\right)+\sum_{\rho_1}\left(\frac{1}{s-\rho_1}+\frac{1}{\rho_1}\right),
\end{align}
where $-a_n$ for $n\geq1$ and $\rho_1:=\beta_1+i\gamma_1$ denote the real and non-real zeros of $\zeta'(s)$, respectively.
It is well known (see for example \cite[Theorem 9]{levmont}) that  $-a_n$ is the only zero of $\zeta'(s)$ in $(-2n-2,-2n)$ and that there are no other zeros in the half-plane $\sigma\leq0$.

Taking the real parts on both sides of \eqref{pfdf}, we obtain
$$
{\rm Re}\frac{\zeta''}{\zeta'}(s)=\frac{\zeta''}{\zeta'}(0)-2+\frac{2(1-\sigma)}{|s-1|^2}+\sum_{\rho_1}{\rm Re}\left(\frac{1}{s-\rho_1}+\frac{1}{\rho_1}\right)+\sum_{n\geq1}\left(\frac{\sigma+a_n}{|s+a_n|^2}-\frac{1}{a_n}\right),
$$
where the terms with $\rho_1$ and $\overline{\rho_1}$ are grouped together (observe that $\zeta'(\overline{\rho_1})=0$). 
Since $\beta_1>0$ unconditionally, we have for any sufficiently large $t>0$ that
\begin{align*}
{\rm Re}\frac{\zeta''}{\zeta'}(s)
&=\frac{\zeta''}{\zeta'}(0)-2+\sum_{\rho_1}\frac{1}{\rho_1}+\sum_{\rho_1}\frac{\sigma-\beta_1}{(\sigma-\beta_1)^2+(t-\gamma_1)^2}-\sum_{n\geq1}\frac{a_n\sigma+|s|^2}{a_n|s+a_n|^2}\\
&\leq\frac{\zeta''}{\zeta'}(0)-2+\sum_{\rho_1}\frac{1}{\rho_1}-\sum_{n\geq1}\frac{\sigma}{|s+a_n|^2}-\sum_{n\leq |s|}\frac{|s|^2}{a_n|s+a_n|^2}.
\end{align*}
The first four terms on the right-hand side of the above relation are $O_\sigma(1)$, while the last sum is of size approximately $-\log t/2$. Hence, the curve of $\zeta(\sigma+it)$ turns in clockwise direction for any sufficiently large $t>0$. If we assume the Riemann hypothesis, then $\beta_1\geq1/2$, as Andreas Speiser had proved \cite{speiser}, and we also obtain Yildirim's result in the left-half of the critical strip.

We have already seen in Theorem \ref{1} that the curvature of $\zeta(\sigma+it)$ alternates its sign infinitely many times as $t\to\infty$ for fixed $1/2<\sigma<1$ by an application of Voronin's theorem. This can be extended to the vertical line $\sigma=1$ by employing another variant of the universality theorem which is implicitly proved by Garunk\v stis {\it et al.} \cite[Theorem 4]{garlamastst}: {\it Let $\sigma_0\in(1/2,1)$, $r>0$, $g:\mathcal{K}\to\mathbb{C}$ continuous, $g(s_0)\neq0$ and analytic for $|s-s_0|<r$. Then, for any $\epsilon\in(0,|g(s_0)|)$, there exist effectively computable positive numbers $T_0$ and $\delta$ such that
\begin{equation}{\label{weakuniv}}
\max_{|s-s_0|\leq\delta r}|\zeta(s+i\tau)-g(s)|<\epsilon
\end{equation}
for any $\tau\in[T,2T]$ and $T\geq T_0$.} The proof does not depend on what $\sigma_0$ may be, which works as a parameter therein, but rather what properties does the curve $\zeta(\sigma_0+it)$ have. In particular, the basic ingredient in the proof of (\ref{weakuniv}) is a weaker version of Voronin's theorem, proved by Voronin \cite{voronin2} himself: {\it If $\sigma_0\in(1/2,1)$, $(a_1,a_2,\dots,a_N)$ is a vector of complex numbers and $\epsilon>0$, then there exists an effectively computable positive number $T_0$ such that
$$
\max_{k\leq N}|\zeta^{(k)}(s_0+i\tau)-a_k|<\epsilon.
$$
for any $\tau\in[T,2T]$ and $T\geq T_0$.} 
However, if we are willing to drop effectivity and just accept the existence of infinitely many $\tau>0$ satisfying the above inequalities, then $\sigma_0=1$ is also admissible as Voronin \cite{voronin1} had shown.
Hence, it can be quickly confirmed that (\ref{weakuniv}) holds also for $\sigma_0=1$ and a divergent sequence of $\tau>0$.

Assume now that $g(s)=e^s$ and $\epsilon>0$ is sufficiently small. Then, it follows by \eqref{weakuniv} that there is some $\delta_+\in(0,1)$ and a divergent sequence of $\tau^+>0$  with
$$
\zeta(s+i\tau^+)=e^s+O(\epsilon),\quad |s-1|\leq\delta_+,
$$
where the big $O$ term is an analytic function in the interior of the disk. 
Therefore,
$$
\zeta'(s+i\tau^+)=e^s+O(\epsilon)=\zeta''(s+i\tau^+),\quad |s-1|\leq\delta_+,
$$
which implies that
$$
{\rm{Re}}\,\frac{\zeta''}{\zeta'}(s+i\tau^+)>0,\quad |s-1|\leq\delta_+.
$$
Thus, the curvature of $\zeta(\sigma+it)$ is positive in disks of center $1+i\tau^+$ and radius $\delta_+$ for a divergent sequence of real numbers $\tau^+>0$. A similar argument with $g(s)=e^{-s}$ shows that the curvature of $\zeta(\sigma+it)$ is negative in disks of center $1+i\tau^-$ and some radius $\delta_-$ for a divergent sequence of real numbers $\tau^->0$. Hence, by taking $\theta=\min\lbrace\delta_+,\delta_-\rbrace$, we conclude that there are infinitely many sign changes of the curvature on the vertical line $1+i\mathbb{R}$ as well as in its immediate right neighbourhood.

\section{Proof of Theorem \ref{3.5}}
{ We begin by setting the necessary theoretical background.}
Assume that $f(s)$ is analytic, not identically zero and almost periodic in the sense of Bohr in a vertical strip $[\alpha,\beta]$.
Then, B\o{}rge Jessen \cite{jessen} has proved that the limit 
$$
\phi_f(\sigma):=\lim_{(\delta-\gamma)\to\infty}\frac{1}{\delta-\gamma}\int_{\gamma}^{\delta}\log |f(\sigma+it)|\mathrm{d}t
$$
exists uniformly in $[\alpha,\beta]$ and represents a continuous and convex function, which is known as the {\it Jensen function} of $f$. As a convex function, $\phi_f(\sigma)$ is almost everywhere  differentiable in the interval $[\alpha,\beta]$ and one can confirm that if $\alpha'\in(\alpha,\beta)$ is such a point, then
\begin{align}\label{limderiv}
\phi_f'(\alpha')=\lim_{(\delta-\gamma)\to\infty}\frac{1}{\delta-\gamma}\int_{\gamma}^{\delta}\mathrm{Re}\frac{f'}{f}(\alpha'+it)\mathrm{d}t.
\end{align}
This result was implicitly shown by Jessen \cite{jessen}, Philip Hartman \cite{hartman}, resp. Jessen and Hans Tornehave \cite{jessentornhave}; we give a proof in the appendix.

In the context of Dirichlet series, Jessen and Tornehave \cite[Theorem 31]{jessentornhave} proved the following:
\begin{Theorem*} [A]
	For an ordinary Dirichlet series $f(s)=\sum_{n\geq n_0}a_{n}n^{-s}$, $a_{n_0}\neq0$, with the uniform convergence abscissa $\alpha$, the Jensen function $\phi_f(\sigma)$ possesses in every half-plane $\sigma>\alpha_1>\alpha$ only a finite number of linearity intervals and a finite number of points of non-differentiability. If $\sigma_0$ denotes the supremum of the real parts of zeros of $f(s)$ (which is always finite), then
	\begin{align}\label{zero-free}
	\phi_f(\sigma)=-\sigma\log n_0+\log|a_{n_0}|,\quad\sigma>\sigma_0.
	\end{align}
	For an arbitrary strip $(\sigma_1,\sigma_2)$ where $\alpha<\sigma_1<\sigma_2<+\infty$, the relative frequency of zeros
	$$
	H_f(\sigma_1,\sigma_2):=\lim_{(\delta-\gamma)\to\infty}\frac{\sharp\left\{\rho:f(\rho)=0,\,\sigma_1<\mathrm{Re}\rho<\sigma_2,\,\gamma<\mathrm{Im}\rho<\delta\right\}}{\delta-\gamma}
	$$
	exists and it is determined by
	$$
	2\pi H_f(\sigma_1,\sigma_2)=\phi_f'(\sigma_2-0)-\phi_f'(\sigma_1+0).
	$$
\end{Theorem*}

This line of research has been pursued shortly afterwards by Vibeke Borchsenius and Jessen  in the case where $f(s)$ can be analytically continued beyond its abscissa of uniform convergence. In particular, if $\phi_f(\sigma)$ and $H_f(\sigma_1,\sigma_2)$ are defined as above with the only difference that $\gamma$ is a sufficiently large but fixed number and only $\delta\to\infty$, then they proved \cite[Theorem 1]{borchjessen} the following:

\begin{Theorem*}[B]
	Let $-\infty\leq\alpha<\alpha_0<\beta_0<\beta\leq+\infty$ and $\gamma\in\mathbb{R}$, and let $f_1(s),f_2(s),\dots$ be a sequence of functions almost periodic in $[\alpha,\beta]$ converging uniformly in $[\alpha_0,\beta_0]$ towards a function $f(s)$.
	Suppose, that none of the functions is identically zero. Suppose further, that $f(s)$ may be continued analytically in the half-strip $\alpha<\sigma<\beta$, $t>\gamma$, and that there is $p>0$ such that
	\[
	\lim_{n\to\infty}\limsup_{\delta\to\infty}\frac{1}{\delta-\gamma}\int_{\gamma}^\delta\int_{\alpha_1}^{\beta_1}|f(\sigma+it)-f_n(\sigma+it)|^p\mathrm{d}\sigma\mathrm{d}t=0,
	\]
	for any reduced strip $\alpha<\alpha_1<\sigma<\beta_1<\beta$.
	Then the Jensen function $\phi_f(\sigma)$ exists in $[\alpha,\beta]$ and it is a continuous and convex function.
	If it is differentiable at $\sigma_1$ and $\sigma_2$ for some $\alpha<\sigma_1<\sigma_2<\beta$, then the relative frequency of zeros $H_f(\sigma_1,\sigma_2)$ exists and it is determined by
	\[
	2\pi H_f(\sigma_1,\sigma_2)=\phi_f'(\sigma_2)-\phi_f'(\sigma_1).
	\]
\end{Theorem*}

\noindent We note that \eqref{limderiv} (with fixed $\gamma$ and $\delta\to\infty$) holds as well in points of the extended interval $[\alpha,\beta]$ where $\phi_f(\sigma)$ is differentiable.
\medskip

Next we apply the aforementioned results in the case of $f=\zeta'$. Firstly, the derivative of the zeta-function can be represented as an absolutely convergent Dirichlet series in the half-plane $\sigma>1$ and it is zero-free in the half-plane $\sigma>E$ for some $E\in[2,3]$ (see \cite[Theorem 11.5 (C)]{tit2}). Hence, in view of \eqref{limderiv} and \eqref{zero-free}  we deduce that
\begin{align}\label{absconv}
\lim_{T\to\infty}\frac{1}{T}\int_{0}^{T}\mathrm{Re}\frac{\zeta''}{\zeta'}(\sigma+it)\mathrm{d}t=\phi_{\zeta'}'(\sigma)=-\log 2,\quad\sigma>E.
\end{align}

For a positive integer $n$, let $\zeta'_n(s):=\sum_{k\leq n}(-\log k)k^{-s}$. Then, from the approximate functional equation for the Riemann zeta-function \cite[Theorem 4.11]{tit2}, i.e.
$$
\zeta(s)=\sum_{k\leq t}\frac{1}{n^s}+O(t^{-\sigma}),\quad t\geq \gamma>0,\quad 0<\sigma_0\leq\sigma\leq3,
$$
it follows that
$$
\lim_{n\to\infty}\limsup_{\delta\to\infty}\frac{1}{\delta-\gamma}\int_{\gamma}^\delta\int_{\alpha_1}^{\beta_1}|\zeta'(\sigma+it)-\zeta'_n(\sigma+it)|^2\mathrm{d}\sigma\mathrm{d}t=0,
$$
for any reduced strip $1/2<\alpha_1<\sigma<\beta_1<3$. Since the functions $\zeta'_n(s)$ for $n\geq1$ are almost periodic in any strip, $\phi_{\zeta'}(\sigma)$ exists in $(1/2,3)$ and, in addition, it is a continuous and convex function. As a matter of fact, it is also infinitely many times differentiable in this interval.

To see this, we observe that
$$
\frac{1}{T}\int_{0}^{T}\log|\zeta'(\sigma+it)|\mathrm{d}t=\frac{1}{T}\int_{0}^{T}\log|\zeta(\sigma+it)|\mathrm{d}t+\frac{1}{T}\int_{0}^{T}\log\left|\frac{\zeta'}{\zeta}(\sigma+it)\right|\mathrm{d}t.
$$
The limit of the first term on the right-hand side is $\phi_\zeta(\sigma)$ and is equal to $0$ for $\sigma>1/2$ (see \cite[Theorem 14, pages 162-163]{borchjessen}). The limit of the second term on the right-hand side has been computed by Charng Rang Guo \cite[Theorem 1.1.3]{guo} and it is equal to
$$
G(\sigma):=\frac{1}{2}\int_{\mathbb{R}}\int_{\mathbb{R}}\widehat{\chi}(x,y;\sigma)\log(x^2+y^2)\mathrm{d}x\mathrm{d}y,
$$
where $\widehat{g}$ denotes the Fourier transform of a function $g$ and
$$
{\chi}(x,y;\sigma):=\prod_p\int_0^1\exp\left(2\pi i\mathrm{Re}((x+iy)a(t,\sigma))\right)\mathrm{d}t,\quad a(t,\sigma):=\log p\sum_{m\geq1}\frac{\exp(2\pi imt)}{p^{m\sigma}}.
$$
Moreover, he showed that \cite[Theorem 1.1.1]{guo} $G(\sigma)$ is infinitely many times differentiable in $(1/2,3]$ with
$$
G^{(k)}(\sigma)=\frac{1}{2}\int_{\mathbb{R}}\int_{\mathbb{R}}\widehat{\frac{\partial^k}{\partial\sigma^k}}\chi(x,y;\sigma)\log(x^2+y^2)\mathrm{d}x\mathrm{d}y,\quad k\geq1.
$$
In other words, Guo found the Jensen function of $\zeta'(s)$ and from his work and our previous discussion follows immediately that
$$
\lim_{T\to\infty}\frac{1}{T}\int_{0}^{T}\mathrm{Re}\frac{\zeta''}{\zeta'}(\sigma+it)\mathrm{d}t=G'(\sigma),\quad\sigma>1/2.
$$
{Moreover, $G'(\sigma)$ is increasing because $G(\sigma)$ is a convex function. 
In combination with \eqref{absconv}, we obtain the theorem.
}

\section{Proof of Theorem \ref{4}}

For fixed $\sigma<1/2$ and $\vert t\vert\geq 2$, Garunk\v{s}tis \& Steuding \cite{garunk} showed that
\begin{equation}\label{lower}
\vert \zeta(\sigma+it)\vert>c\vert t\vert^{1/2-\sigma+\epsilon},
\end{equation}
where $\epsilon>0$ is arbitrary and $c>0$ is an absolute constant depending only on $\epsilon$ and $\sigma$. This lower bound relies on the conditional estimate $\zeta(\sigma+it)\gg \vert t\vert^{-\epsilon}$ for fixed $\sigma>1/2$ (see \cite{tit2}, \S 14.2) in combination with the functional equation. Since the right-hand side of (\ref{lower}) tends with $\vert t\vert$ to infinity, the curve $t\mapsto \zeta(\sigma+it)$ {\it spirals} outside the disk of radius $r$ centered at the origin for all sufficiently large values of $\vert t\vert$. More precisely, if $\vert t\vert\geq r^{2/(1-2\sigma)}$, then, for $r\geq r_0$, the values $\zeta(\sigma+it)$ lie in the complement of the disk. Therefore, it only remains to show that the image of the values $\zeta(\sigma+it)$ for $t$ from a finite interval, $[0,T]$ say, cannot be locally dense. 

Assume that the values $\zeta(\sigma+ it)$ are somewhere dense for $t\in[0, T ]$. Then there exists a small disk $D$ such that the curve $t \to \zeta(\sigma+it)$ visits for $t\in[0, T ]$ every tiny disk $d$ of radius $(3N)^{-1}$ centered at the elements of $\frac{1}{N} \mathbb{Z}[i]\cap D$, where $N$ is a large integer. Since the tiny disks $d$ have distance $\gg N^{-1}$ one from another and there are $\gg N^2$ many of them in $D$, the curve $t \to\zeta(\sigma+it)$ has length $\gg N^2\cdot N^{-1}=N$ which tends with $N$ to infinity. This contradiction implies that the curve is nowhere dense.

\section{Appendix}

If
$$
\phi_f(\sigma;\gamma,\delta):=\frac{1}{\delta-\gamma}\int_{\gamma}^{\delta}\log |f(\sigma+it)|\mathrm{d}t,
$$
then it follows from Jessen and Tornhave \cite[page 187]{jessentornhave} that
\begin{eqnarray}\label{fourineq}
\phi_f'(\sigma-0)& \leq& \mathop{\underline{\lim}}_{(\delta-\gamma)\to\infty}\phi_f'(\sigma-0;\gamma,\delta)\nonumber\\
&\leq & \left\{\begin{array}{l}\mathop{\overline{\lim}}_{(\delta-\gamma)\to\infty}\phi_f'(\sigma-0;\gamma,\delta)\\
\mathop{\underline{\lim}}_{(\delta-\gamma)\to\infty}\phi_f'(\sigma+0;\gamma,\delta)
\end{array}\right\}\\
&\leq&\mathop{\overline{\lim}}_{(\delta-\gamma)\to\infty}\phi_f'(\sigma+0;\gamma,\delta)\leq\phi_f'(\sigma+0).\nonumber
\end{eqnarray}
Therefore, if $\phi_f(\sigma)$ is differentiable at a point $\sigma=\alpha'$, then all of the above four limits are equal. 
In particular, we have
\begin{align}\label{deriv}
\phi_f'(\alpha'-0;\gamma,\delta)=\phi_f'(\alpha'+0;\gamma,\delta)+o(1)=\phi_f'(\alpha')+o(1)
\end{align}
as $\delta-\gamma\to\infty$.

Now let $\gamma<\delta$ be a pair of real numbers. If $f(s)$ has no zeros in the vertical segment $[\alpha'+i\gamma,\alpha'+i\delta]$, then it is analytic in a neighbourhood of this segment and, consequently, $\phi_f(\sigma;\gamma,\delta)$ is differentiable at $\sigma=\alpha'$ with
\begin{align*}
\phi_f'(\alpha';\gamma,\delta)	
&=\frac{\mathrm{d}}{\mathrm{d}\sigma}\left[\frac{1}{\delta-\gamma}\int_{\gamma}^{\delta}\log|f(\alpha'+it)|\mathrm{d}t\right]_{\sigma=\alpha'}\\
&=\frac{1}{\delta-\gamma}\int_{\gamma}^{\delta}\left.\frac{\partial\log|f(\sigma+it)|}{\partial\sigma}\right|_{\sigma=\alpha'}\mathrm{d}t=\frac{1}{\delta-\gamma}\int_{\gamma}^{\delta}\mathrm{Re}\frac{f'(\alpha'+it)}{f(\alpha'+it)}\mathrm{d}t.
\end{align*}
Next we assume that $f(s)$ has only one zero in the vertical segment $[\alpha'+i\gamma,\alpha'+i\delta]$ of order $m$, $\alpha'+it_0$ say. Then there is some $\varepsilon>0$ such that
$$
f(s)=(s-\alpha'-it_0)^mg(s),\quad|s-\alpha'-it_0|\leq\varepsilon,
$$
where $g(s)$ is analytic and zero-free in the disk $|s-\alpha'-it_0|<2\varepsilon$. Therefore,
\begin{align*}
\frac{\mathrm{d}}{\mathrm{d}\sigma}&\left[\int_{t_0-\varepsilon}^{t_0+\varepsilon}\log|f(\sigma+it)|\mathrm{d}t\right]_{\sigma=\alpha'\pm0}=\\
&=\frac{\mathrm{d}}{\mathrm{d}\sigma}\left[\int_{t_0-\varepsilon}^{t_0+\varepsilon}m\log|\sigma-\alpha'+i(t-t_0)|+\log|g(\sigma+it)|\mathrm{d}t\right]_{\sigma=\alpha'\pm0}\\
&=\lim_{\sigma\to\alpha'^\pm}\int_{t_0-\varepsilon}^{t_0+\varepsilon}\frac{m(\sigma-\alpha')\mathrm{d}t}{(\sigma-\alpha')^2+(t-t_0)^2}+\int_{t_0-\varepsilon}^{t_0+\varepsilon}\mathrm{Re}\frac{g'(\alpha'+it)}{g(\alpha'+it)}\mathrm{d}t\\
&=\pm m\pi+\int_{t_0-\varepsilon}^{t_0+\varepsilon}\mathrm{Re}\frac{f'(\alpha'+it)}{f(\alpha'+it)}\mathrm{d}t.
\end{align*}
It thus follows that 
$$
\phi_f'(\alpha'\pm0;\gamma,\delta)=\frac{1}{\delta-\gamma}\int_{\gamma}^{\delta}\mathrm{Re}\frac{f'(\alpha'+it)}{f(\alpha'+it)}\mathrm{d}t\pm\frac{m\pi}{\delta-\gamma}.
$$
{ If $f$ had more than one different zeros in $[\alpha'+i\gamma,\alpha'+i\delta]$, then $m$ would represent the sum of all such zeros counted with multiplicity.}
In view of \eqref{deriv} we deduce first that $m\pi=o(\delta-\gamma)$, $\delta-\gamma\to\infty$, and then we obtain the desired relation \eqref{limderiv}.

The same reasoning applies in the strip of analytic continuation, where relation \eqref{fourineq} still holds assuming the conditions of Theorem (B) (see Borchsenius and Jessen \cite[page 122]{borchjessen}).

\bigskip

\noindent
Athanasios Sourmelidis\\
Institute of Analysis and Number Theory, TU Graz\\
Steyrergasse 30, 8010 Graz, Austria\\
sourmelidis@math.tugraz.at
\medskip

\noindent
J\"orn Steuding\\
Department of Mathematics, W\"urzburg University\\
Emil Fischer-Str. 40, 97\,074 W\"urzburg, Germany\\
steuding@mathematik.uni-wuerzburg.de


\begin{thebibliography}{10}

\bibitem{ander}{\sc J. Andersson}, Lavrent$'$ev's approximation theorem with nonvanishing polynomials and universality of zeta-functions, in: {\it New Directions in Value-Distribution theory of Zeta and L-Functions}, Shaker Verlag, Aachen, 2009, 7–10.

\bibitem{bagchi}{\sc B. Bagchi}, {\it The Statistical Behaviour and Universality Properties of the Riemann Zeta-Function and Other Allied Dirichlet Series}, Ph.D. Thesis, Calcutta, Indian Statistical Institute, 1981. 

\bibitem{bohr}{\sc H. Bohr, R. Courant}, Neue Anwendungen der Theorie der Diophantischen Approximation auf die Riemannsche Zetafunktion, {\it J. reine angew. Math.} {\bf 144} (1914), 249-274.

\bibitem{borchjessen}{\sc V. Borchsenius, B. Jessen}, Mean motions and values of the Riemann zeta function.
{\it Acta Math.} {\bf 80} (1948), 97–166. 

\bibitem{carmo}{\sc M.P. do Carmo}, {\it Differential Geometry of Curves and Surfaces}, Prentice-Hall, 1976. 

\bibitem{edwards}{\sc H.M. Edwards}, {\it Riemann's Zeta Function}, Academic Press, New York-London, 1974.

\bibitem{fujii}{\sc A. Fujii}, On a conjecture of Shanks, {\it Proc. Jpn. Acad.} {\bf 70} (1994), 109-114. 

\bibitem{garlamastst}{\sc R. Garunk\v stis, A. Laurin\v cikas, K. Matsumoto, J. Steuding, R. Steuding}, Effective uniform approximation by the Riemann zeta-function, {\it Publ. Mat.} {\bf 54} (2010), 209–219. 

\bibitem{garunk}{\sc R. Garunk\v{s}tis, J. Steuding}, On the roots of the equation $\zeta(s)=\alpha$, {\it Abh. Math. Seminar Univ. Hamburg} {\bf 84} (2014), 1-15. 

\bibitem{gone}{\sc S.M. Gonek}, {\it Analytic Properties of Zeta and $L$-Functions}, Ph.D. Thesis, University of Michigan, 1979.

\bibitem{gonek}{\sc S.M. Gonek, H.L. Montgomery}, Spirals of the Zeta Function I, in: {\it Analytic Number Theory. In Honour of Helmut Maier's 60th Birthday}, {\sc C. Pomerance, M.Th. Rassias} (eds.), Springer, 2015, 127-131. 

\bibitem{guo}{\sc C.R. Guo}, On the zeros of the derivative of the Riemann zeta function, {\it Proc. London Math. Soc.} (3) {\bf 72} (1996), no. 1, 28–62. 

\bibitem{hartman}{\sc P. Hartman}, Mean motions and almost periodic functions, {\it Amer. Math. Soc.} {\bf 46}, (1939). 66–81. 

\bibitem{jessen}{\sc B. Jessen}, \"Uber die Nullstellen einer analytischen fastperiodischen Funktion. Eine Verallgemeinerung der Jensenschen Formel. (German) {\it Math. Ann.} {\bf108} (1933), no. 1, 485–516. 

\bibitem{jessentornhave}{\sc B. Jessen, H. Tornehave}, Mean motions and zeros of almost periodic functions.
{\it Acta Math.} {\bf 77} (1945), 137–279. 

\bibitem{lapidus}{\sc H. Herichi, M.L. Lapidus}, {\it Quantized Number Theory, Fractal Strings, and Riemann Hypothesis}, World Scientific, 2021.

\bibitem{kalpokas}{\sc J. Kalpokas, J. Steuding}, On the value-distribution of the Riemann zeta-function on the critical line, {\it Mosc. J. Comb. Number Theory} {\bf 1} (2011), 26-42.

\bibitem{karatsuba}{\sc A.A. Karatsuba, S.M. Voronin}, {\it The Riemann Zeta-Function}, de Gruyter, 1992.

\bibitem{lapi}{\sc M.L. Lapidus}, Towards quantized number theory: spectral operators and an asymmetric criterion for the Riemann hypothesis, {\it Philos. Trans. Royal Soc., Ser. A} {\bf 373} (2015), 24 pp.

\bibitem{levmont}{\sc N. Levinson,  H.L. Montgomery}, Zeros of the derivatives of the Riemann zeta-function,
{\it Acta Math.} {\bf 133} (1974), 49–65. 

\bibitem{mergelyan}{\sc S.N. Mergelyan}, Uniform approximations to functions of a complex variable, {\it Usp. Mat. Nauk} {\bf 7} (1952), 31-122 (Russian); {\it Am. Math. Soc.} {\bf 101} (1954). 

\bibitem{netto}{\sc E. Netto}, Beitrag zur Mannigfaltigkeitslehre, {\it J. reine angew. Math.} {\bf 86} (1879), 263-268.

\bibitem{sagan}{\sc H. Sagan}, {\it Space-Filling Curves}, Springer 1994. 

\bibitem{shanks}{\sc D. Shanks}, Review of Haselgrove, {\it Math. Comput.} {\bf 15} (1961), 84-86. 

\bibitem{speiser}{\sc A. Speiser}, Geometrisches zur Riemannschen Zetafunktion, {\it Math. Ann.} {\bf 110} (1935), 514–521.

\bibitem{19}{\sc A. Sourmelidis, J. Steuding}, An atlas for all plane curves, submitted to {\it EMS Magazine}

\bibitem{steu}{\sc J. Steuding}, {\it Value-Distribution of $L$-Functions}, Lecture Notes in Mathematics {\bf 1877}, Springer 2007.

\bibitem{tit2}{\sc E. C. Titchmarsh}, {\it The Theory of the Riemann Zeta-Function}, Oxford University Press, 2nd ed. revised by D.R. Heath-Brown, 1986.

\bibitem{voronin1}{\sc S.M. Voronin}, The distribution of the nonzero values of the Riemann $\zeta$-function (Russian), {\it Collection of articles. Part II, Trudy Mat. Inst. Steklov} {\bf 128} (1972), 131–150.

\bibitem{voronin}{\sc S.M. Voronin}, Theorem on the `universality' of the Riemann zeta-function, {\it Izv. Akad. Nauk SSSR, Ser. Matem.} {\bf 39} (1975), 475-486 (Russian); {\it Math. USSR Izv.} {\bf 9} (1975), 443--445.

\bibitem{voronin2}{\sc S.M. Voronin}, $\Omega$-theorems of the theory of the Riemann zeta-function (Russian), {\it Izv. Akad. Nauk SSSR Ser. Mat.} {\bf 52} (1988), 424–436, 448; {\it Math. USSR-Izv.} {\bf 32} (1989), 429–442. 

\bibitem{wegert0}{\sc E. Wegert}, {\it Visual Complex Functions}, Birkh\"auser, 2012. 

\bibitem{wegert}{\sc E. Wegert, G. Semmler}, Phase Plots of Complex Functions: a Journey in Illustration, {\it Notices A.M.S.} {\bf 58} (2011), 768-780.

\bibitem{yildirim}{\sc C. Y\i ld\i r\i m}, A note on $\zeta''(s)$ and $\zeta'''(s)$, {\it Proc. Am. Math. Soc.} {\bf 124} (1996), 2311-2314. 

\end{thebibliography}
\end{document}